\documentclass{article}
\usepackage{amsthm}
\usepackage{amssymb, amsmath, bbm}
\newtheorem{Theorem}{Theorem}
\newtheorem{Lemma}{Lemma}
\newtheorem{Proposition}{Proposition}
\newcommand{\R}{{\mathbb{R}}}
\newcommand{\p}{{\mathcal{P}}}

\newcommand{\A}{{\mathcal{A}}}
\numberwithin{equation}{section}
\numberwithin{Theorem}{section}
\numberwithin{Lemma}{section}
\numberwithin{Proposition}{section}
\usepackage{array}
\begin{document}
\title{A Positive Semidefinite Approximation of the Symmetric Traveling Salesman Polytope}

\author{Ellen Veomett}
\date{October, 2006}

\maketitle

\begin{abstract} For a convex body $B$ in a vector space $V$, we construct its approximation $P_k, k=1, 2, \dotsc$
using  an intersection of a cone of positive semidefinite quadratic forms with an affine subspace.   We show that $P_k$
is contained in $B$ for each $k$.  When $B$ is the Symmetric Traveling Salesman Polytope on  $n$ cities $T_n$, we show
that the scaling of $P_k$ by $\frac{n}{k}+ O\bigl(\frac{1}{n}\bigr)$ contains $T_n$ for  $k \leq \lfloor \frac{n}{2}
\rfloor$.  Membership for $P_k$ is computable in time polynomial in $n$ (of degree linear in $k$).   

We discuss facets of $T_n$ that lie on the boundary of $P_k$.  We introduce a new measure on each facet defining
inequality for $T_n$ in terms of the eigenvalues of a quadratic form.  Using these eigenvalues of facets, we show 
that the scaling of $P_1$ by $O(\sqrt{n})$ has all of the facets of $T_n$ defined by the subtour elimination constraints
either in its interior or lying on its boundary.
\end{abstract}

\section{Introduction and Results}

For many interesting convex bodies $X$ in a vector space $V$, given a point $x \in V$, the question ``is $x$ in $X$?'' is
difficult to answer.
  This fact has generated work in the direction of finding another set $Y$ which is ``close'' to $X$ in some way for which
the  membership question is ``easy'' to answer.  Sherali and Adams \cite{Sherali-Adams}, Lov\'{a}sz and Schrijver 
\cite{Lovasz-Schrijver},  and Lasserre  \cite{Lasserre:2} have constructed approximating sets in the case where the body
to  be approximated is a 0-1 polytope.  In each of these instances, the authors constructed successive relaxations of a 0-1 polytope,
 such that in the $n$th step, the 0-1 polytope is achieved: $P = K^n \subset K^{n-1} \subset \dotsb \subset K^1 \subset
K$.   Metric properties of these approximating sets are not known.  For specifics, as well as a comparison of the methods, see \cite{Laurent}.  

In the following we construct successive relaxations $P_1, P_2, \dotsc$ of an arbitrary convex body $X$, each of which is
contained in
$X$.  If $X \subset \R^n$ is a 0-1 polytope, then we also obtain $P_n = X$ (for details see Section \ref{0-1}).   We
explore in particular the case where $X$ is the Symmetric Traveling Salesman Polytope, where we estimate the closeness of the approximation metrically.

\subsection{The Symmetric Traveling Salesman Polytope}\label{STSPdescription}

The Symmetric Traveling Salesman Polytope (STSP) can be described as follows: recall that a Hamiltonian cycle in the complete graph on $n$ vertices $K_n$ is a cycle which visits every vertex exactly once.  To each Hamiltonian cycle in $K_n$, we can associate its incidence matrix $A = (a_{ij})$ where
\begin{equation*}
a_{ij}= \begin{cases}
1 & \text{if the cycle contains edge } \{i,j\} \\
0 & \text{if the cycle does not contain edge } \{i,j\}
\end{cases}
\end{equation*}
The polytope is called symmetric because there is a similar notion in the case of a digraph (a graph where the edges have an orientation), which is the Asymmetric Traveling Salesman polytope.

Note that each matrix corresponding to a Hamiltonian  cycle is a symmetric 0-1 matrix in $\R^{n^2}$ with 0s on the
diagonal.  Given a particular matrix corresponding to a  Hamiltonian cycle, any other such matrix can be obtained from
it by simultaneously permuting rows and columns (this  corresponds to permuting the labels on the vertices of the
graph).  The Symmetric Traveling Salesman Polytope is typically described as the convex  hull of all adjacency matrices
corresponding
to Hamiltonian cycles in $K_n$.  However, for our purposes, instead of using the entire adjacency matrix for each
Hamiltonian cycle, we will only use the
``upper half'' of each such matrix.  Thus, we consider the STSP to be the convex hull of the
$\frac{n(n-1)}{2}$-dimensional vectors indexed by pairs $\{i,j\}$ ($i,j \in \{1, 2, \dotsc, n\}$) where entry
$x_{ij}=1$ if the associated Hamiltonian cycle contains edge $\{i,j\}$ and $x_{ij}=0$ if it does not contain edge
$\{i,j\}$.

Thus, each vertex of the STSP is  a vector which corresponds to a cycle. 
To each cycle we can
associate a permutation of the numbers $\{1, 2, \dotsc, n\}$  beginning with the number 1 where the permutations $(1,
m_2, m_3, \dotsc, m_n)$ and $(1, m_n, m_{n-1}, \dotsc, m_2)$  are identified.  We will use the descriptions of the
vertices as vectors, cycles, and permutations interchangeably.   Using the permutation description of a Hamiltonian
cycle, it is not hard to see that there are $\frac{(n-1)!}{2}$ different Hamiltonian cycles in $K_n$.
The STSP has been studied widely, though a complete description  of it via linear inequalities is not known (and in
some sense, cannot be known unless NP=coNP, see \cite{PapK}).   It is clearly not full dimensional in $\R^{n(n-1)/2}$;
for example, for each point in the STSP, the sum of the entries is $n$.   It is not hard to show that its dimension is
$\frac{n(n-3)}{2}$.   For more information on the STSP and the associated  Traveling Salesman Problem, see, for
example, Chapter 58 of \cite{Schrijver:book}.  Linear optimization over the  STSP and the membership question for the
STSP are known to be NP-hard.

\subsection{Semidefinite Construction}\label{cone1}

The following observation of A. Barvinok \cite{B:personal} gives the construction with which we will work.  Let $V$ be a real vector space and let $X \subset V$ be a finite (though possibly very large) set.  Let $V^*$ denote the dual of $V$.  Recall that the polar dual of $X$ is the set
\begin{equation*}
X^\circ = \{f:X \to \R : \quad f \text{ is linear}, \quad f(x) \leq 1 \text{ for all } x \in X\} \subset V^*
\end{equation*}
We note that by ``$f$ linear'' we mean that $f$ is the restriction to $X$ of a linear function on $V$.   We view $X^\circ$ as living in
 the space $\R^X$ of all functions from $X$ to $\R$.    If the convex hull of $X$ does not contain the origin in its relative
interior,  $X^\circ$ is not bounded.  Indeed,  one can find a linear function $f \in V^*$ not identically 0 on all of $X$
which  separates the origin and $X$ ($f(x)\leq 0 \text{ for all }
 x \in X$) so that $\alpha f \in X^\circ$ for all nonnegative $\alpha$.   We thus consider the polar of $X$ in  its
affine span with the center of polarity being the barycenter of $X$:
\begin{equation*} 
A = \{f:X \to \R :\;\;  f \text{ is affine}, \;\; f(x) \leq 1 \text{ for all } x \in X, \;\; \frac{1}{|X|} \sum_{x \in X} f(x) = 0 \}
\end{equation*}
Again, we note that by ``$f$ affine'' we mean that $f$ is the restriction to $X$ of an affine function on $V$.  Then for convenience, we flip $A \mapsto -A$ and then shift $f \mapsto f+1$ so that we obtain the following description of the dual:
\begin{equation*}
Q = \{ f:X \to \R :  \;\;  f \text{ is affine}, \;\; f(x) \geq 0 \text{ for all } x \in X, \;\; \frac{1}{|X|} \sum_{x \in X} f(x) = 1 \}
\end{equation*}
The set with which we will work is $Q$.  

We note that any convex body $B$ can be written as the polar dual to some other convex body $B'$ which is in the same affine span, with the center of polarity the barycenter of $B'$.  Since $B'$ can be arbitrarily closely approximated by a convex hull of finitely many points, $B$ is arbitrarily close to some $Q$ as defined above.

Fix a positive integer $k$ and let $\p_k(V)$ be the space of all polynomials of degree at most $k$ on $V$.  To any function $f:X \to \R$ we can associate the quadratic form 
\begin{equation*}
q_f: \p_k(V) \to \R
\end{equation*}
defined by
\begin{equation*}
q_f(h) = \frac{1}{|X|}\sum_{x \in X}f(x)h^2(x) \quad \text{for} \quad h \in \p_k(V)
\end{equation*}
Clearly, if $f(x) \geq 0$ for each $x \in X$, then $q_f$ is a positive semidefinite quadratic form on $\p_k(V)$.

Note that, as $f$ ranges over affine functions on $X$ with average value 1, the form $q_f$ ranges over  an affine
subspace in the space of quadratic forms on $\p_k(V)$.  Let us define $\A_k$ to be this affine subspace,  and define  $W_k$ to
be  the cone of positive semidefinite quadratic forms $q: \p_k(V) \to \R$.  We define $P_k = \{f: q_f \in \A_k \cap W_k\}$.  Then we can see that
\begin{equation*}
Q \subset P_k
\end{equation*}
which leads us to ask
\begin{itemize}
\item How close is $P_k$ to $Q$?
\end{itemize}

\subsection{The Case of a 0-1 Polytope}\label{0-1}
We note that if $X \subset \R^n $ consists of 0-1 vectors, then $P_n = Q$.  Indeed, let $f \in P_n$ so that $q_f$ is a positive semidefinite quadratic form and $f$ corresponds to an affine function with average value 1 on $X$.  Let us fix any $y \in X$.  Let $I \subset \{1, 2, \dotsc, n\}$ consist of the indices of the entries of $y$ which are 0, and $J \subset \{1, 2, \dotsc, n\}$ be the indices of the entries of $y$ which are 1.  Then we can see that the degree $n$ polynomial
\begin{equation*}
p_y(x) = \prod_{i \in I}(1-x_i)\prod_{j \in J}x_j
\end{equation*}
has value 1 on $y$ and 0 on any other vector in $X$.  Thus, we have
\begin{equation*}
0 \leq q_f(p_y) = \frac{1}{|X|}\sum_{x \in X} f(x)p_y^2(x) = \frac{f(y)}{|X|}
\end{equation*}
Since $y$ was arbitrary, we see that $f(x) \geq 0$ for each $x \in X$, so that $f \in Q$, giving us $P_n \subset Q$.  Since we already had $Q \subset P_n$, we see that indeed $Q = P_n$.

\subsection{The Case of the STSP}

From this point on, we fix $X$ to be the set of vectors corresponding to Hamiltonian cycles in $K_n$ as described in
section \ref{STSPdescription}.  Thus, here our vector space is $\R^{n(n-1)/2}$.  Recall that each vector in
$\R^{n(n-1)/2}$ is indexed by unordered pairs $\{i,j\}$ where $i,j \in \{1, 2, \dotsc, n\}$.  For $x \in
\R^{n(n-1)/2}$, we use $x_{ij}$ and $x_{ji}$ interchangeably to denote the entry of $x$ corresponding to the pair
$\{i,j\}$.  
The barycenter of the STSP is the vector $Z = (z_{ij})$ where 
\begin{equation*}
z_{ij} = \begin{cases} \frac{2}{n-1} & \text{if } i \not= j \\ 0 & \text{if } i=j
\end{cases}
\end{equation*}
and the average value of any affine function on $X$ is simply its value on $Z$.  Note that the all ones function: $\mathbbm{1}(x) =
1$   for all $x \in X$ can be written as a linear function, since in the affine span of $X$ it corresponds to the inner
product with the vector $(\frac{1}{n}, \frac{1}{n}, \dotsc, \frac{1}{n})$.  Thus, if $f$ is a linear function and $a \in \R$, then in the affine span of $X$ the affine function $f(x)+a$ is equal to the linear function $f(x)+a\mathbbm{1}(x)$.  Hence we can see that in the set 
\begin{equation*} 
A = \{f:X \to \R :\;\;  f \text{ is affine}, \;\; f(x) \leq 1 \text{ for all } x \in X, \;\; \frac{1}{|X|} \sum_{x \in X} f(x) = 0 \}
\end{equation*}
as defined in Section \ref{cone1}, we can actually have each $f$ being \textit{linear}.  Flipping $A'\mapsto -A'$ doesn't destroy linearity of the functions, nor does shifting $f \mapsto f+1$.  Hence, in this
case, we have
\begin{equation*}
Q = \{ f:X \to \R :  \;\;  f \text{ is linear}, \;\; f(x) \geq 0 \text{ for all } x \in X, \;\; \frac{1}{|X|} \sum_{x \in X} f(x) = 1 \}.
\end{equation*}

We define $P_k$ just as in Section \ref{cone1}, so that $P_k$ is the set of all linear functions $f$ with average value 1 on $X$ whose corresponding quadratic form $q_f$ is positive semidefinite:
\begin{align*}
 0 \leq q_f(p) = \frac{1}{|X|}\sum_{x \in X} f(x)p^2(x)  \quad \quad & \begin{split} & \text{for polynomials } p(x)
\text{ on }
\R^{n(n-1)/2} \\
& \text{of degree } \leq k  \end{split}
\end{align*}
Note that the function $\mathbbm{1}$ is the center of $Q$.  We have the following:

\begin{Theorem}\label{semidefinite1}
For any $n\geq 9$ and any $k=1, 2, \dotsc, \lfloor \frac{n}{2} \rfloor$, there exists a constant $a_k = \frac{n}{k} + \alpha_k$ where $|\alpha_k| \leq \frac{c}{n}$ ($c$ can be taken to be  10 for any $n$ and $k$) such that 
\begin{equation*}
Q -\mathbbm{1}\; \subset \; P_k -\mathbbm{1} \; \subset \; a_k(Q-\mathbbm{1})
\end{equation*}
\end{Theorem}

Recall that by our definition of $Q$, the traveling salesman polytope $T_n$ is the set of points in the affine span of $X$ such that $f(x) -1 \geq -1$ for all $f \in Q$.  Thus, defining $P_k^\circ$ to be the set of points in the affine span of $X$ such that $g(x)-1 \geq -1$ for all $g \in P_k$, and denoting the barycenter of $T_n$ by $Z$, we have

\begin{Theorem}\label{semidefinite2}
Let $T_n$ be the symmetric traveling salesman polytope.  For any $n \geq 9$ and any $k =1, 2, \dotsc, \lfloor \frac{n}{2} \rfloor$, there exists a constant $a_k = \frac{n}{k} + \alpha_k$ where $|\alpha_k| \leq \frac{c}{n}$ ($c$ can be taken to be  10 for any $n$ and $k$) such that 
\begin{equation*}
P_k ^\circ -Z \; \subset \; T_n -Z \; \subset \; a_k(P_k^\circ -Z)
\end{equation*}
\end{Theorem}

Note that the approximation gives us an upper bound on how far $P_k$ is from $T_n$.   We also note that approximating
the STSP with respect to its center gives us a new measure of approximation for the Traveling Salesman Problem:
approximation with respect to the average value.
Specifically, suppose that $f$ is a linear objective function on the STSP.
  Then using Theorem \ref{semidefinite2}, we
can bound the difference between the optimal value  and average value of $f$ on  $T_n$ based on the difference between
the optimal value of $f$ on $P_k^\circ$ and the average value of $f$ on $T_n$.

Let $C$ be the cone of positive semidefinite quadratic forms on a vector space $V$.  Then membership  in $C$ is
decidable in time of order $\dim^3(V)$ (see, for example Chapter 1 of \cite{GLS}).  In the case where $V =
\p_k(\R^{n(n-1)/2})$, $\dim(V) = {n(n-1)/2 +k \choose k}$ so that membership in $K_k$ is decidable in time of order $n^{6k}$.

The remainder of the paper is structured as follows: in Section \ref{boundary} we discuss facets  of the STSP which we
know lie on the boundary of $P_k$, and in Section \ref{UsingEigenvalues} we use eigenvalues to analyze the first
approximation $P_1$ with respect to the subtour elimination constraints.  In section
\ref{computations} we prove the bounds in Theorem \ref{semidefinite1} and in section \ref{EigenvalueCalculation} we
prove the bounds discussed in section \ref{UsingEigenvalues}.

\section{Facets on the Boundary}\label{boundary}

Although there is no known complete description of the Symmetric Traveling Salesman Polytope as a system of linear inequalities, many facets are known (see, for example, chapter 58 of \cite{Schrijver:book}).  Some well-known facet defining inequalities are the following:

\begin{align}
 0  \leq x_{ij} & \leq 1   &\text{ for each } i,j \label{facet1} \\
 \sum_{\substack{j \in U \\ i \in V-U}} x_{ij} & \geq 2  &\text{ for each } U \subset V \text{ with } \emptyset \not= U \not= V \label{facet2} \\
 \sum_{\substack{j \in U  \\ i \in V-U \\ \{i,j\} \not\in F}} x_{ij} -\sum_{\{k,\ell\} \in F} x_{k\ell}  & \geq 1-|F|      \begin{split} & \text{for } U \subset V, \quad F \text{ matching}, \\
& |F| \geq 3 \text{ odd}, \text{ each edge of } F  \\  & \text{ having one endpoint in } U \end{split} \label{facet3}
\end{align}
Inequalities \eqref{facet2} are known as the subtour elimination constraints, and inequalities \eqref{facet3} are
known as the 2-matching constraints.

Any
facet of the STSP can defined by some linear inequality $f(x)
\geq 0$ which is unique up to a scaling. 
If we
scale so that the average value on $X$ is 1, then we know that the scaled function must be in $P_k$.   A natural
question to ask would be: which (if any) of the linear functions defining a facet for the STSP lie on the \textit{boundary} of $P_k$?  

Let $h_{ij}$ be the linear function such that $h_{ij}(x) \geq 0$ corresponds to the right hand  side of inequality
\eqref{facet1} for edge $\{i,j\}$.  Consider the degree 1 polynomial $p_{ij} = x_{ij}$.   Then $h_{ij}(x)$ is 0
whenever $x$ contains the edge $\{i,j\}$, and $p_{ij}(x)$ is 0 whenever $x$ does not contain the edge $\{i,j\}$.  Thus, we have
\begin{equation*}
q_{h_{ij}}(p_{ij}) = \frac{1}{|X|}\sum_{x \in X}h_{ij}(x)(p_{ij}(x))^2 = 0
\end{equation*}
so that $h_{ij}$ lies on the boundary of $P_k$ for any $k \geq 1$.  If we let $h_{ij}'$ be the polynomial 
corresponding to the left hand side of inequality \eqref{facet1} for the edge $\{i,j\}$, then we can again easily see
\begin{equation*}
q_{h_{ij}'}(1-p_{ij}) = \frac{1}{|X|}\sum_{x \in X}h_{ij}'(x)(1-p_{ij}(x))^2 = 0
\end{equation*}
so that $h_{ij}'$ lies on the boundary of $P_k$ for any $k \geq 1$.

Now suppose that $h_U$ is a linear function such that $h_U(x) \geq 0$  corresponds to the facet \eqref{facet2} for some
$U \subset V, \emptyset \not= U \not= V$.  Let $|U| = m$ and $U = \{\ell_1, \ell_2, \ldots, \ell_m\}$.  Consider the
degree $m-1$ polynomial 
\begin{equation*}
p_U = x_{\ell_1\ell_2}x_{\ell_2\ell_3}\dotsm x_{\ell_{m-1}\ell_{m}}
\end{equation*}
Note that $h_U(x) =0$ whenever $U$ has two outgoing edges, and  $p_U(x) \not= 0$ only if $x$ contains the path $\ell_1,
\ell_2, \dotsc, \ell_m$, which implies that there are two edges going out of $U$.  Thus, we have
\begin{equation*}
q_{h_U}(p_U) = \frac{1}{|X|}\sum_{x \in X}h_U(x)p_U^2(x) = 0
\end{equation*}
so that $h_U$ is on the boundary of $P_k$ if $k \geq m-1$

In general, suppose we have  a linear function $f\geq 0$ that defines a facet.   Then if we can construct a polynomial
of degree $\leq k$ such that any cycle for which that polynomial is nonzero  must be on the facet defined by $f$, $q_f$
is on the boundary of $P_k$.  For example, consider $U$ and $F$ as in  \eqref{facet3}, say  $U =\{\ell_1, \ell_2, \dotsc,
\ell_m\}$
 and $F = \{\{\ell_1, n_1\}, \{\ell_2, n_2\}, \dotsc, \{\ell_{2s+1}, n_{2s+1}\}\}$ ($3 \leq 2s+1 \leq m$).  Let $h_{U,F}$ 
be
a linear function such that $h_{U,F} \geq 0$ defines the facet in \eqref{facet3} corresponding to $U$ and $F$.   Define
the degree $s+m$ polynomial $p_{U,F}(x)$ as follows:
\begin{equation*}
p_{U,F}(x) =  \prod_{i=1}^{2s+1}x_{\ell_in_i} \prod_{j=1}^{s}x_{\ell_{2j-1}\ell_{2j}} \prod_{k=2s+1}^{m-1}x_{\ell_k\ell_{k+1}}
\end{equation*}
Then whenever $p_{U,F}(x) \not=0$, $x$ contains the paths $n_1\ell_1\ell_2n_2, \; n_3\ell_3\ell_4n_4,$ \linebreak $
\dotsc,\; n_{2s-1}\ell_{2s-1}\ell_{2s}n_{2s}, \; n_{2s+1}\ell_{2s+1}\ell_{2s+2}\ell_{2s+3}\dotsm \ell_m$.  This implies
that 
\begin{equation*}
 \sum_{\substack{j \in U  \\ i \in V-U \\ \{i,j\} \not\in F}} x_{ij} = 1
 \end{equation*}
 and
 \begin{equation*}
\sum_{\{\ell_i,n_i\} \in F} x_{\ell_i, n_i} = 2s+1 = |F|
\end{equation*}
so that $h_{U,F}(x) = 0$.  Thus, 
\begin{equation*}
q_{h_{U,F}}(p_{U,F}) = \frac{1}{|X|}\sum_{x \in X}h_{U,F}(x)p_{U,F}^2(x) = 0
\end{equation*}
so that $h_{U,F}$ is on the boundary of $P_k$ if $k \geq s+m$.

\section{Approximation Appraisal via Eigenvalues}\label{UsingEigenvalues}

Recall from Theorem \ref{semidefinite1} that when we considered the vector space of polynomials of degree no more than 1
on $\R^{n(n-1)/2}$, if $Q$ is the dual to the TSP and $P_1$ is the first approximating set, our bounds were
\begin{equation*}
Q -\mathbbm{1} \subset P_1 -\mathbbm{1} \subset n(Q-\mathbbm{1})
\end{equation*}
A reasonable question to ask is: can we find a scaling factor $a\ll n$ such that we still have $P_1 - \mathbbm{1} \subset
a(Q-\mathbbm{1})$?  
Using eigenvalues we will
see that, if the TSP had no more facets than the subtour elimination constraints (defined in section \ref{boundary}),
the answer to the above question would be yes with $a = \sqrt{n}$.

We will use the same notation for the subtour elimination constraints that we used in section \ref{boundary}, denoting
by $h_U$ the function corresponding to subset $U \subset V$, where $m = |U| \leq \frac{n}{2}$. 
Since we are restricting ourselves here to the first approximating set $P_1$,
we consider the quadratic form $q_{h_U}$ only on the space of polynomials on $\R^{n(n-1)/2}$ of degree no more than 1.

We firstly note that for any $x \in X$ and any $1 \leq i \leq n$, the polynomial 
\begin{equation*}
2-\sum_{\substack{1 \leq j \leq n \\ \j \not= i}} x_{ij}
\end{equation*}
has value 0 (this corresponds to the fact that, in any Hamiltonian  cycle, each vertex has exactly two edges incident
to it).  Thus,  for any $f:X \to \R$, we will consider the vector  corresponding to that polynomial to be an
eigenvector with eigenvalue 0 for $q_f$.  There are $n$ such linearly independent eigenvectors; one for each vertex.

Let $M$ be the set of functions corresponding to the subtour elimination constraints, and suppose we only  wanted to
have the functions $a(M-\mathbbm{1})$ be outside $P_1-\mathbbm{1}$ or on its boundary for some $a>0$.   (Note that this
is \textit{necessary} in order to have $ P_1 -\mathbbm{1} \subset a(Q-\mathbbm{1})$, as each function  in $M$ is on the
boundary of $Q$).  This is equivalent to requiring that, for each function $h_U$, the matrix corresponding to the 
quadratic form $q_f$ with $f = ah_U+(1-a)\mathbbm{1}$ has a 0 or negative eigenvalue (beyond the $n$ 0 eigenvalues  which
each quadratic form is known to have).

We have the following:

\begin{Theorem}\label{eigenvalues}
Consider  the quadratic form $q_f$ with $f = ah_U+(1-a)\mathbbm{1}$ defined on the space of 
linear functions on $\R^{n(n-1)/2}$.   
For
each $U\subset V$  $|U|=m, 3 \leq m \leq \frac{n}{2}$, if $n = |V| \geq 6$, then $q_{ah_U+(1-a)\mathbbm{1}}$ has the
following eigenvalues:

\begin{tabular}{ll}
\textbf{Multiplicity} & \textbf{Eigenvalue} \\
  n &   0 \\
  $\frac{m(m-3)}{2}$ &    $a\frac{2(m-2)}{(n-2)(m-1)} + (1-a)\frac{2}{n-1}$ \\
  $\frac{(n-m)(n-m-3)}{2}$ &    $a\frac{2(n-m-2)}{(n-2)(n-m-1)} + (1-a)\frac{2}{n-1}$ \\
  $(n-m-1)(m-1)$ &   $a\frac{2(mn^2-nm^2-n^2+4n-3mn+3m^2-4)}{(n-2)(n-3)(mn-m^2-n+1)} + (1-a)\frac{2}{n-1}$ \\
  $m-1$ &   $a\frac{2(m-2)}{(n-3)(m-1)} + (1-a)\frac{2}{n-1}$ \\
  $n-m-1$&   $a\frac{2(n-m-2)}{(n-3)(n-m-1)} + (1-a)\frac{2}{n-1}$ \\
  1 &   $\alpha+\sqrt{\beta}$  \\
  1 &   $\alpha-\sqrt{\beta}$
\end{tabular} 

where $\alpha$ and $\beta$ are rational functions in $a, m$, and $n$ such that if $a = \sqrt{n}$, then for all $m \leq
\frac{n}{2}$ we have $\alpha - \sqrt{\beta} \leq 0$.
\end{Theorem}

We note that  if we want to force one of the first five nonzero eigenvalues to be nonpositive for any $m \leq \frac{n}{2}$,
we would need $a = O(n)$.
If one calculates the exact expressions for $\alpha$ and $\beta$ (see section \ref{EigenvalueCalculation}), one sees
that for a smaller 
$m$, a smaller $a$ is required to force the last eigenvalue to be nonpositive. 
This corresponds to our intuition
that the facets defined by $h_U$ for small $|U|=m$ are big and deep facets, because they contain many vertices.  
Thus, intuitively it should be easier to find polynomials $p$ for small $m=|U|$ for which $\frac{1}{|X|}\sum_{x \in X}h_U(x)p^2(x)
\approx 0$, implying that that the functions $h_U$ are closer to the boundary of $P_1$.

As noted above, the fact that one of the eigenvalues for $q_{ah_U+(1-a)\mathbbm{1}}$ is nonpositive for $a = \sqrt{n}$
implies that
if we scale each of the subtour elimination constraints by $\sqrt{n}$, the function lies outside or on the boundary of
$P_1$.  This fact   has several possible
implications.  It could mean that the bound of $n$ from Theorem \ref{semidefinite1} is closer to the real bound, but that
we need to look for facets beyond the
subtour elimination constraints in order to see that the bound of $n$ is necessary.  Or it could mean that the bound of $\sqrt{n}$ is closer to the real bound, and we have yet to find a way to prove this.

We note that there is a polynomial time separation algorithm for the subtour elimination constraints.  Indeed, an $x
\in \R^{n(n-1)/2}$ satisfies the subtour elimination constraints if and only if the minimum cut for the complete graph
with capacities corresponding to the entries of $x$ is at least 2.  The author is unaware of any known lift
constraints whose description is polynomial in $n$ and whose projection achieves the subtour elimination constraints.

\section{Proofs of Metric Bounds}\label{computations}

Recall that from Section \ref{cone1}, we already have $Q \subset P_k$.  Thus, the only question we must address is: how far is $P_k$ from $Q$?  In other words, given a function $f \in P_k$, that is, a function defining a positive semidefinite quadratic form via
\begin{equation*}
q_f(h) = \frac{1}{|X|}\sum_{x \in X}f(x)h^2(x) \quad \text{for} \quad h \in \p_k(\R^{n(n-1)/2})
\end{equation*}
where $f$ is a linear function with average value 1 on $X$, how negative can the values of $f$ on $X$ be?  The following Lemma gives us a bound:

\begin{Lemma}\label{bound}
Fix $y \in X$ and $f \in P_k$ so that $f$ is a linear function with average value 1 on $X$ and $q_f$ is positive semidefinite.  Suppose that we find polynomials $p_{1}, \ldots, p_{m}$ of degree $k$ such that $p_{i}$ takes on only values 0 or 1 and there exist positive constants $b_k < c_k$ such that for any $i,j \in \{ 1, 2, \dotsc n\}, i \not= j$,
\begin{equation*}
\sum_{x\in X:\{i,j\} \in x} \sum_{i=1}^{m} p_{i}(x) = \begin{cases}
b_k & \text{ if } \{i,j\} \not \in y \\
c_k & \text{ if } \{i,j\} \in y
\end{cases}
\end{equation*}
Then 
\begin{equation*}
-\frac{b_k(n-1)}{2(c_k-b_k)} \leq f(y)
\end{equation*}
\end{Lemma}

\begin{proof}
Let $f \in P_k$ so that $f:X \to \R$ is linear function with average value 1 on $X$.  Let us fix $y \in X$.  For each $i < j$, we define the 
vector $e_{ij} = (\epsilon_{st}) \in \R^{n(n-1)/2}$ as follows:
\begin{equation*}
\epsilon_{st}=\begin{cases} 1 & \text{ if } \{s,t\}=\{i,j\} \\
0 & \text{ otherwise } \end{cases}
\end{equation*}
Note that each $x \in X$ can be written as a sum of the vectors $e_{ij}$ ($i<j$) for which $\{i,j\}$ is an edge in $x$. 
Each $e_{ij}$ will appear in exactly $(n-2)!$ different $x \in X$.  Thus, the fact that $f$ has average 1 on $X$ tells us:

\begin{align*}
1 &= \frac{2}{(n-1)!} \sum_{x \in X} f(x) \\
&= \frac{2}{(n-1)!} \sum_{x \in X}f\left(\sum_{\{i,j\} \in x, i<j} e_{ij} \right)\\ 
&=\frac{2}{(n-1)!} \sum_{i<j} f(e_{ij})(n-2)! \\
&= \sum_{i<j} \frac{2}{n-1}f(e_{ij})
\end{align*}
which gives us
\begin{equation}\label{1}
\frac{n-1}{2}-f(y) = \sum_{\{i,j\} \not\in y}f(e_{ij})
\end{equation}
for any particular $y \in X$.

Since $f \in P_k$, the form $q_f$ is positive semidefinite, so for any polynomial $p(x)$ we can write the inequality
\begin{equation*}
0 \leq q_f(p) = \frac{2}{(n-1)!}\sum_{x \in X} f(x)p^2(x)
\end{equation*}
which implies
\begin{equation*}
0 \leq \sum_{x \in X}f(x)p^2(x) =  \sum_{\substack{ x \in X \\ i<j \\ \{i,j\} \in x}} f(e_{ij})p^2(x)
\end{equation*}

Now assuming we have $p_i$ as stated in the Lemma, for each $k$ we find that 
\begin{equation*}
 0 \leq\sum_{\substack{ x \in X \\ i<j \\ \{i,j\} \in x}} f(e_{ij})p_{\ell}^2(x) = \sum_{\substack{ x \in X \\ i<j \\ \{i,j\} \in x}} f(e_{ij})p_{\ell}(x) \quad \text{ for } \quad \ell = 1, 2, \ldots m
\end{equation*}
so that using \eqref{1} we have
\begin{align*}
 0 \leq\sum_{\substack{ x \in X \\ i<j \\ \{i,j\} \in x}} \sum_{\ell=1}^{m} f(e_{ij})p_{\ell}(x)  &=  c_k\sum_{\{i,j\} \in y} f(e_{ij}) + b_k\sum_{\{i,j\} \not\in y}f(e_{ij}) \\
 & = c_kf(y) +b_k(\frac{n-1}{2} -f(y))
 \end{align*}
 which then implies
 \begin{equation*}
 -\frac{b_k(n-1)}{2(c_k-b_k)} \leq f(y)
\end{equation*}

\end{proof}

We note that Lemma \ref{bound} only gives a bound on how negative a function $f \in P_k$ can be, if we can find polynomials $p_i$ satisfying the assumptions.  It may be that, in fact, $f$ is entirely nonnegative.  Picking a particular set of polynomials, we will prove the following:

\begin{Proposition}\label{calculations} Let us fix $y \in X$ and $f \in P_k$.  If $n$ is even, then
\begin{equation*}
 -\frac{n}{k}+1-\frac{n(k-1)}{k(n^2-kn-3n+k+3)} \leq f(y)
 \end{equation*}
If $n$ is odd then
 \begin{equation*}
 -\frac{n}{k}+1- \frac{n(k-1)}{k(n^2-nk-4n+4+2k)}  \leq f(y)
 \end{equation*}
\end{Proposition}
To complete the calculations required for Proposition \ref{calculations}, we need a Lemma:

\begin{Lemma}\label{paths}
Let $(k_1, k_2, \dotsc, k_m)$ be a partition of $k$ ($k+m \leq n$) and $K_n $ the complete graph on $n$ vertices.  Let $p_1, p_2, \dotsc, p_m$ be disjoint paths in $K_n$ of length $k_1, \dotsc, k_m$ respectively.  Then the number of Hamiltonian cycles in $K_n$ containing all of paths $p_1, \dotsc, p_m$ is:
\begin{equation*}
2^{m-1}(n-k-1)!
\end{equation*}
\end{Lemma}

\begin{proof}

Note that the restriction $k+m \leq n$ assures us that it is possible to find disjoint paths in $K_n$ of lengths $k_1, \dots, k_m$.  Any cycle containing the paths $p_1, \dotsc, p_m$ can be written uniquely as a sequence of numbers, beginning with path $p_1$ in a particular orientation.  Thinking of the remaining paths as blocks with 2 orientations and the remaining numbers as blocks with a single orientation, we find that there are $2^{m-1}(n-k-1)!$ ways of ordering and orienting the remaining blocks.  Each of these orders and orientations corresponds uniquely to a Hamiltonian cycle containing paths $p_1, \dotsc, p_m$.
\end{proof}

\begin{proof}[Proof of Proposition \ref{calculations}.]
We will use Lemma \ref{bound}.  First we need to describe the polynomials which we will use.   Note that in the Hamiltonian cycle $y$, depending on whether $n$ is either even or odd, there are either two or $n$ different subsets of $\lfloor \frac{n}{2} \rfloor$ disjoint edges in $y$.  We will call such a subset of edges an ``EO subset of  $y$'' (EO for ``every other'').  For each EO subset $\Gamma$ of $y$ and each $I \subset \Gamma$ of cardinality $k$, we define:
\begin{displaymath}
p_{I,\Gamma} = \prod_{\{i,j\} \in I} x_{i,j}
\end{displaymath}
In words, $p_{I,\Gamma}$ is the monomial corresponding to $k$ disjoint edges which are a subset of some EO subset of $y$. 
(We note that this is why we need to restrict $k \leq \lfloor \frac{n}{2} \rfloor$).  

Note that each $p_{I,\Gamma}$ takes
on only values 0 or 1. 
In order to use Lemma
\ref{bound}, we will need to calculate
\begin{equation}
\label{3}
\sum_{I,\Gamma} \sum_{x:\{i,j\} \in x} p_{I,\Gamma}(x)
\end{equation}
where in the first sum  $\Gamma$ runs over all EO subsets of $y$, and $I$ runs over all $k$-element subsets.  We note that these polynomials were chosen with Lemma \ref{bound} in mind; namely so that for each edge $\{i,j\}$, \eqref{3} has only two different values: one value if $\{i,j\} \not\in y$ and another value if $\{i,j\} \in y$.

Suppose that $n$ is even.  Then $y$ has two EO subsets, $\Gamma_1$ and $\Gamma_2$.  Note that when we calculate (\ref{3}), we are simply counting the number of Hamiltonian cycles containing both some $I \subset \Gamma_\ell$ of size $k$ and the edge $\{i,j\}$.  Note that in each of $\Gamma_1$ and $\Gamma_2$, for each $i \in \{1, 2, \ldots, n\}$, there is exactly one edge which contains $i$.  If $\{i,j\} \not\in y$, the edge which contains $i$ and the edge which contains $j$ are distinct.  If $\{i,j\} \in y$ then $\{i,j\}$ is in one of $\Gamma_1$ or $\Gamma_2$.  In the other, the edge which contains $i$ and the edge which contains $j$ are distinct.

Let us pick some edge $\{i,j\} \not\in y$.  Then for \textit{each} of the EO subsets $\Gamma_1$ and $\Gamma_2$  there are ${n/2 -2 \choose k-2}$ subsets $I$ of size $k$ containing the edge which contains $i$ and the edge which contains $j$.  For such subsets $I$, $I \cup \{i,j\}$  consists of $k-1$ distinct paths, $k-2$ of which are of length 1, and 1 of which is of length 3.  

There are $2{n/2-2 \choose k-1}$ subsets $I$ of size $k$ containing exactly one of the edge which contains $i$ or the edge which contains $j$.  For such subsets $I$, $I \cup \{i,j\}$ consists of $k$ distinct paths, $k-1$ of which are of length 1, 1 of which is of length 2.

Lastly, there are  ${n/2-2 \choose k}$ subsets $I$ of size $k$ containing neither the edge which contains $i$ nor the edge which contains $j$.  For such subsets $I$, $I \cup \{i,j\}$ consists of $k+1$ distinct paths, each of length 1.  Thus, from Lemma \ref{paths}, we can see that if $\{i,j\} \not \in y$ then we can calculate \eqref{3} (which we denote $f_1(n,k)$) to be 

\begin{gather}
f_1(n,k) =\sum_{I,\Gamma} \sum_{x:\{i,j\} \in x} p_{I,\Gamma}(x) \notag \\
=  2 \left[ {\frac{n}{2} -2 \choose k-2} 2^{k-2} (n-k-2)! + 2{\frac{n}{2}-2 \choose k-1}2^{k-1}(n-k-2)! \right. \notag \\
  \left. + {\frac{n}{2}-2 \choose k}2^k(n-k-2)!\right] \label{even1}
\end{gather}

Recall that if $\{i,j\} \in y$,  exactly one of $\Gamma_1$ or $\Gamma_2$ contains the edge $\{i,j\}$, say $\Gamma_1$ does.  Then $\Gamma_2$ contains 1 edge which contains $i$, and a disjoint edge which contains $j$.  By arguments similar to those above, and again using Lemma \ref{paths},  we can see that if $\{i,j\}  \in y$ then we can calculate \eqref{3} (which we denote $f_2(n,k)$) to be 

\begin{gather}
f_2(n,k) = \sum_{I,\Gamma} \sum_{x:\{i,j\} \in x} p_{I,\Gamma}(x) \notag \\
={\frac{n}{2}-1 \choose k-1}2^{k-1}(n-k-1)! +{\frac{n}{2}-1 \choose k}2^k(n-k-2)! \notag \\
 + {\frac{n}{2} -2 \choose k-2} 2^{k-2} (n-k-2)! + 2{\frac{n}{2}-2 \choose k-1}2^{k-1}(n-k-2)! \notag \\
 + {\frac{n}{2}-2 \choose k}2^k(n-k-2)! \label{even2}
\end{gather}

Thus, using these calculations and Lemma \ref{bound}, we see that if $n$ is even and $f \in P_k$ then
\begin{gather*}
-\frac{(n-1)}{2}\frac{f_1(n,k)}{f_2(n,k)-f_1(n,k)} \\
= -\frac{n}{k}+1-\frac{n(k-1)}{k(n^2-kn-3n+k+3)} \leq f(y)
\end{gather*}

Now suppose that $n$ is odd.  Then $y$ has $n$ EO subsets, $\Gamma_1, \ldots, \Gamma_n$, where $\Gamma_i$ does not have an edge coming from vertex $i$.  

Note that for each $i \in \{1, 2, \ldots, n\}$, and each $\Gamma_\ell, \ell \not= i$, there is exactly one edge which contains $i$.  If $\{i,j\} \not\in y$ and $i,j\not= \ell$, then in $\Gamma_\ell$ the edge which contains $i$ and the edge which contains $j$ are distinct.  

If $\{i,j\} \in y$  then $\{i,j\}$ is in $ \frac{n-1}{2}$ of the $\Gamma_j$s.  In $ \frac{n-1}{2}  -1$ of the $\Gamma_j$s the edge containing $i$ and the edge containing $j$ are distinct.  And in $\Gamma_j$, there is only an edge containing $i$, in $\Gamma_i$ there is only an edge containing $j$. 

Let us pick some edge $\{i,j\} \not\in y$.  Then for the EO subsets $\Gamma_\ell, \ell \not= i,j$  there are ${(n-1)/2 -2 \choose k-2}$ subsets $I$ of size $k$ containing the edge which contains $i$ and the edge which contains $j$.  For such subsets $I$, $I \cup \{i,j\}$ consists of $k-1$ disjoint paths, $k-2$ of which are of length 1, 1 of which is of length 3.   There are $2{(n-1)/2-2 \choose k-1}$ subsets $I$ of size $k$ containing exactly one of the edge which contains $i$ or the edge which contains $j$.  For such subsets $I$, $I \cup \{i,j\}$ consists of $k$ disjoint paths, $k-1$ of which are of length 1, 1 of which is of length 2.  And there are ${(n-1)/2-2 \choose k}$ subsets $I$ of size $k$ containing neither the edge which contains $i$ nor the edge which contains $j$.  For such subsets $I$, $I \cup \{i,j\}$ consists of $k+1$ disjoint paths, each of length 1. 

In $\Gamma_i$, there are ${(n-1)/2 -1\choose k-1}$ subsets $I$ of size $k$ containing the edge which contains $j$ ($I \cup \{i,j\}$ consisting of $k-1$ paths of length 1, 1 path of length 2), and ${(n-1)/2 -1\choose k}$ subsets $I$ of size $k$ not containing the edge which contains $j$ ($I \cup \{i,j\}$ consisting of $k+1$ paths of length 1).  Similarly, in $\Gamma_j$, there are ${(n-1)/2 -1\choose k-1}$ subsets $I$ of size $k$ containing the edge which contains $i$ ($I \cup \{i,j\}$ consisting of $k-1$ paths of length 1, 1 path of length 2), and ${(n-1)/2 -1\choose k}$ subsets $I$ of size $k$ not containing the edge which contains $i$ ($I \cup \{i,j\}$ consisting of $k+1$ paths of length 1).

Recall that in calculating (\ref{3}), we are simply counting the number of Hamiltonian cycles containing both some $I$ of size $k$ and the edge $\{i,j\}$.  Thus, from Lemma \ref{paths}, we can see that if $\{i,j\} \not \in y$ then we can calculate \eqref{3} (which we denote $g_1(n,k)$) to be 
\begin{gather}
g_1(n,k) = \sum_{I,\Gamma} \sum_{x:\{i,j\} \in x} p_{I,\Gamma}(x) \notag \\
= (n-2) \left[ {\frac{n-1}{2} -2 \choose k-2} 2^{k-2} (n-k-2)! + 2{\frac{n-1}{2}-2 \choose k-1}2^{k-1}(n-k-2)! \right. \notag \\
\left. + {\frac{n-1}{2}-2 \choose k}2^k(n-k-2)!\right] \notag \\
 + 2\left[ {\frac{n-1}{2} -1\choose k-1} 2^{k-1}(n-k-2)! + {\frac{n-1}{2} -1\choose k}2^k(n-k-2)! \right] \label{odd1}
\end{gather}

Recall that if $\{i,j\} \in y$,  $ \frac{n-1}{2} $ of the $\Gamma_\ell$s contain the edge $\{i,j\}$, $\frac{n-1}{2} -1$ of the $\Gamma_\ell$s have the edge containing $i$ and the edge containing $j$ being distinct, $\Gamma_j$ does not have an edge which contains $j$ and $\Gamma_i$ does not have an edge which contains $i$.  By arguments similar to those above, and again using Lemma \ref{paths}, we find that for $\{i,j\} \in y$ we can calculate \eqref{3} (which we denote $g_2(n,k)$) to be 
\begin{gather}
g_2(n,k) = \sum_{I,\Gamma} \sum_{x:\{i,j\} \in x} p_{I,\Gamma}(x) \notag \\
= \frac{n-1}{2}  \left[{ \frac{n-1}{2}-1 \choose k-1}2^{k-1}(n-k-1)! +{\frac{n-1}{2}-1 \choose k}2^k(n-k-2)! \right] \notag \\
 +\left(  \frac{n-1}{2}-1\right) \left[ {\frac{n-1}{2} -2 \choose k-2} 2^{k-2} (n-k-2)! + 2{\frac{n-1}{2}-2 \choose k-1}2^{k-1}(n-k-2)! \right. \notag \\
 \left. + {\frac{n-1}{2}-2 \choose k}2^k(n-k-2)! \right]  \notag \\
 + 2 \left[ {\frac{n-1}{2} -1\choose k-1} 2^{k-1}(n-k-2)! + {\frac{n-1}{2} -1\choose k}2^k(n-k-2)! \right] \label{odd2}
\end{gather}

Thus, using these calculations and Lemma \ref{bound}, we see that if $n$ is odd and $f \in P_k$ then 
\begin{gather*}
-\frac{(n-1)}{2} \frac{g_1(n,k)}{g_2(n,k)-g_1(n,k)}\\
= -\frac{n}{k}+1- \frac{n(k-1)}{k(n^2-nk-4n+4+2k)}  \leq f(y)
\end{gather*}

\end{proof}

Now we can prove Theorem \ref{semidefinite1}:

\begin{proof}[Proof of Thm \ref{semidefinite1}]
Recall that we assume $n \geq 9$ and $\lfloor \frac{n}{2} \rfloor \geq k$.  Note that both 
\begin{equation*}
\frac{n(k-1)}{k(n^2-kn-3n+k+3)} \quad \quad \text { and } \quad \quad \frac{n(k-1)}{k(n^2-nk-4n+4+2k)}\end{equation*}
 are bounded above in absolute value by $\frac{c}{n}$ for an absolute constant $c$ (which can, for example, be 10).  
Thus, from Proposition \ref{calculations} we know that there exists $a_k = \frac{n}{k} + \alpha_k$ with $|\alpha_k| \leq
\frac{c}{n}$  for an absolute constant $c$ such that, if $f \in P_k$,  for each $y \in X$, $-a_k+1 \leq f(y)$.   This
implies that $(f+(a_k-1)\mathbbm{1})(y) \geq 0$ for all $y \in X$.  It is clear that $f+(a_k-1)\mathbbm{1}$ has 
average value $a_k$ on $X$ (recall that $f$ has average value 1 on $X$).  It is also clear that $f+(a_k-1)\mathbbm{1}$ 
is a linear function on $X$ (recall that $f$ is linear; the function $\mathbbm{1}$ corresponds to the inner product with
the  vector $(\frac{1}{n}, \frac{1}{n}, \dotsc, \frac{1}{n})$).  Thus, we have $f+(a_k-1)\mathbbm{1} \in a_kQ$.   Thus,
we have 
\begin{equation*}
Q -\mathbbm{1} \; \subset \; P_k -\mathbbm{1} \; \subset \; a_k(Q-q_\mathbbm{1})
\end{equation*}
\end{proof}

\section{Eigenvalues of Facets}\label{EigenvalueCalculation}

Here we work out the calculations needed for Theorem \ref{eigenvalues}.  Recall that for a function $f: X
\to\R$, we considered the quadratic form $q_f$ defined on the vector space $\p_k(\R^{n(n-1)/2})$ of polynomials of
degree no more than $k$ on $\R^{n(n-1)/2}$.  Before we restrict ourselves to $k=1$ (so that we consider $q_f$ to be
only on the space of linear functions on $\R^{n(n-1)/2}$), we prove a theorem for general $k$:

\begin{Theorem}\label{trace}
Let $f:X \to \R$ be a function with average value 1 on $X$: $\frac{1}{|X|}\sum_{x \in X}f(x) = 1$.  Consider the
quadratic form $q_f$ on \linebreak $\p_k(\R^{n(n-1)/2})$.  Then

\begin{equation*}
\text{trace}(q_f) =  {n+k \choose k}
\end{equation*}
\end{Theorem}

\begin{proof}
We will use the matrix $A_f$ associated to this quadratic form with respect to the orthonormal basis of monomials $\{1,
x_{12}, x_{13}, \dotsc, x_{ij}, \dotsc, x_{12}x_{13},
\dotsc\}$.  \linebreak 
The diagonal entries of $A_f$ will be 
\begin{equation*}
\frac{1}{|X|}\sum_{x
\in X}f(x)\Big(\prod_{\alpha
\in I} x_{\alpha}\Big)^2(x)
\end{equation*}
 where $I$
is some multiset of edges of size no more than $k$.  (Thus the monomial $x_{12}^2x_{34}$ corresponds to the multiset $\{\{1,2\},\{1,2\},\{3,4\}\}$ of size 3).  Note that 
\begin{equation*}
\Big(\prod_{\alpha \in I} x_{\alpha}\Big)^2(x) = \begin{cases} 0 & \text{ if } x \text{ does not contain all edges in } I
\\
1 &
\text{ if } x \text{ contains all edges in } I 
\end{cases}
\end{equation*}
Thus, the trace of $A_f$ is
\begin{equation*}
\sum_{I} \frac{1}{|X|}\sum_{\substack{x \in X \\ I \subset x}}f(x)
\end{equation*}
where the first sum is over all multisets $I$ of edges such that $|I| \leq k$.  For each Hamiltonian cycle $x$, there are
exactly ${n+k \choose k}$ multisets of edges of size no more than $k$ which are in $x$.  Thus, we can see that 
\begin{align*}
\sum_{I} \frac{1}{|X|}\sum_{\substack{x \in X \\ I \subset x}}f(x) &= \frac{1}{|X|}\sum_{x \in X}{n+k \choose k} f(x) \\
&={n+k \choose k}
\end{align*}
\end{proof}

Note that the above proof actually shows that 
\begin{equation*}
\text{Trace}(A_f) = {n+k \choose k}\Big(\frac{1}{|X|}\sum_{x \in X}f(x)\Big).
\end{equation*}

Now we shall restrict ourselves to considering quadratic forms on the space of linear functions on $\R^{n(n-1)/2}$. 
First we
inspect the entries of the matrices corresponding to these quadratic forms with respect to the
orthonormal basis of monomials $\{1, x_{12}, x_{13}, \dotsc \}$.  
Let $f$ be a real valued function on $X$ and consider a linear function on $\R^{n(n-1)/2}$:
\begin{equation*}
p(x) = \alpha_0+ \sum_{1 \leq i < j \leq n} \alpha_{ij} x_{ij} 
\end{equation*}
We note that
\begin{align*}
q_{f}(p) &= \frac{1}{|X|} \sum_{x \in X}f(x)p^2(x) \\
&= \frac{1}{|X|} \sum_{x \in X} f(x)\bigg(\alpha_0+ \sum_{1 \leq i < j \leq n} \alpha_{ij} x_{ij}\bigg)^2 \\
&= \frac{1}{|X|} \sum_{x \in X} f(x) \bigg(\sum_{i<j, p<q} \alpha_{ij}\alpha_{pq}x_{ij}x_{pq} + \sum_{i<j} \alpha_0\alpha_{ij}x_{ij} + \alpha_0^2\bigg) \\
&=\sum_{i<j, p<q} \alpha_{ij}\alpha_{pq} \bigg(\frac{1}{|X|} \sum_{x \in X} f(x) x_{ij}x_{pq} \bigg) + \sum_{i<j} \alpha_0\alpha_{ij} \bigg( \frac{1}{|X|} \sum_{x \in X} f(x)x_{ij}\bigg) \\
&+\alpha_0^2\bigg(\frac{1}{|X|} \sum_{x \in X} f(x)\bigg)
\end{align*}
For each $x \in X$, each term $x_{ij}$ is either 1 or 0, depending on whether or not $x$ contains edge $\{i,j\}$.  Thus, we can see
that  the entries of the matrix $A_f$ corresponding to $q_f$ are $\frac{1}{|X|}\sum_{I \subset x} f(x)$ where $I$ is some
0, 1, or 2 element subset of the edges in $K_n$.

Recall from section \ref{boundary} that for $U \subset V$ with $m = |U| \leq \frac{n}{2}$, the subtour elimination
constraint $h_U$ can be defined as follows: 
\begin{equation*}
h_U(x) = c\bigg(\sum_{\substack{u \in U \\ v \in V-U}} x_{uv} - 2\bigg)
\end{equation*}
where $c$ is chosen so that the average value of $h_U(x)$ on $X$ is 1.   

We will prove a couple of Lemmas before proving Theorem \ref{eigenvalues}

\begin{Lemma}\label{h_U}
Let $A_U$ be the matrix corresponding to the quadratic form $q_{h_U}$ acting on the vector space of linear functions
on $\R^{n(n-1)/2}$.  
Then for
each $U\subset V$  $|U|=m, 3 \leq m \leq \frac{n}{2}$, if $n = |V| \geq 6$,  $A_U$ has the following eigenvectors and
eigenvalues:

\begin{tabular}{p{1.5in}p{2.5in}}
\textbf{Eigenvalue} & \textbf{Eigenvector} \\ \hline
 0& $2-\sum_{\substack{1 \leq j \leq n \\ \j \not= i}} x_{ij}$, any $1 \leq i \leq n$\\ \hline
  $\frac{2(m-2)}{(n-2)(m-1)} $& $x_{ij}-x_{jf}+x_{fg}-x_{g i}$, distinct $i,j,f,g \in U$   \\ \hline
 $\frac{2(n-m-2)}{(n-2)(n-m-1)}$& $x_{pq}-x_{qr}+x_{rs}-x_{sp}$, distinct $p,q,r,s \in V-U$   
\\ \hline
$\frac{2(m(n-3)(n-m)-(n-2)^2)}{(n-2)(n-3)(m-1)(n-m-1)}$ & $x_{ip}-x_{iq}+x_{jq}-x_{jp}$, distinct $i,j \in U$,
distinct $q,p \in V-U$   \\ \hline
$\frac{2(m-2)}{(n-3)(m-1)}$& $\sum_{\substack{\ell \in U \\ \ell \not= i,j}}  \bigg( \frac{n-m}{m-2}x_{i\ell}
-\frac{n-m}{m-2} x_{j\ell}\bigg) + \sum_{t \in V-U} \bigg(-x_{it}+x_{jt}\bigg)$, distinct $i,j \in U$   \\ \hline
$\frac{2(n-m-2)}{(n-3)(n-m-1)}$ &    $\sum_{\substack{t \in V- U \\ t \not= p,q}} \bigg( \frac{m}{n-m-2} x_{pt} 
-\frac{m}{n-m-2} x_{qt}\bigg) + \sum_{\ell \in U} \bigg(-x_{p\ell}+x_{q\ell}\bigg)$, distinct $p,q \in V-U$
\end{tabular}   
\end{Lemma}

\begin{Lemma}\label{ones}
Let $A_\mathbbm{1}$ be the matrix corresponding to the quadratic form $q_{\mathbbm{1}}$ acting on the vector space of
linear functions on $\R^{n(n-1)/2}$ where $\mathbbm{1}$ is the ``all ones function'': $\mathbbm{1}(x) = 1$ for all $x
\in X$.
Then  $A_\mathbbm{1}$ has the following
eigenvectors and eigenvalues:

\begin{tabular}{p{1.4in}p{2.5in}}
\textbf{Eigenvalue} & \textbf{Eigenvector} \\ \hline
 0& $2-\sum_{\substack{1 \leq j \leq n \\ \j \not= i}} x_{ij}$, any $1 \leq i \leq n$\\ \hline
$\frac{2}{n-1} $ & $x_{\alpha \beta}-x_{\beta \gamma}+x_{\gamma \delta} - x_{\delta \alpha}$ distinct $\alpha, \beta, \gamma, \delta
\in V $
\end{tabular}   
\end{Lemma}

\begin{proof}[Proof of Lemma \ref{h_U}]  Firstly we note that we have already shown that there are $n$ linearly
independent eigenvectors with
eigenvalue 0.  Before we prove that the 5 remaining eigenvectors are, in fact, eigenvectors, we need to be able to
calculate the entries in $A_U$

Recall that the entries of $A_U$ are 
\begin{align}
\frac{1}{|X|}\sum_{\substack{x \in X \\ I \subset x}} h_U(x) &= \frac{1}{|X|}\sum_{\substack{x \in X \\ I \subset x}} c\bigg(\sum_{\substack{u \in U \\ v \in V-U}} x_{uv} - 2\bigg)  \notag \\
&= \frac{c}{|X|}\bigg(\sum_{\substack{u \in U \\ v \in V-U}}\sum_{\substack{x \in X \\ I \subset x}}x_{uv} - \sum_{\substack{x \in X \\ I \subset x}}2 \bigg) \label{entryA_U}
\end{align}
where $I$ is some 0,1, or 2 element subset of the edges.   Note that the sum in \eqref{entryA_U} depends only on whether the endpoints of the edges in $I$ are in $U$ or $V-U$, and how the edges overlap; it does not depend on the labels of the vertices.  In other words, the sum in \eqref{entryA_U} is invariant under the action of $S_m \times S_{n-m}$
on the edges in $I$, where $S_m$ permutes only the indices corresponding to vertices in $U$ and $S_{n-m}$ permutes
indices corresponding to the vertices in $V-U$.

Let us calculate the value of $c$.  Recall that it is chosen so that $h_U$ has an average value of 1 on $X$.  Thus, using Lemma \ref{paths}, we calculate:
\begin{align*}
& \frac{1}{|X|}\sum_{x \in X}\bigg(\sum_{\substack{u \in U \\ v \in V-U}} x_{uv} - 2 \bigg) \\
= & \frac{2}{(n-1)!}\sum_{\substack{u \in U \\ v \in V-U}} \sum_{x \in X} x_{uv}-\frac{2}{(n-1)!}\sum_{x \in X}2 \\
= & \frac{2}{(n-1)!} \sum_{\substack{u \in U \\ v \in V-U}} \#\{\text{Hamiltonian cycles containing } \{u,v\}\} - 2 \\
= & \frac{2}{(n-1)!} \sum_{\substack{u \in U \\ v \in V-U}} (n-2)! -2 \\
= & \frac{2}{(n-1)!}m(n-m)(n-2)!-2 \\
= & \frac{2(m(n-m)+1-n)}{n-1}
\end{align*}
so that 
\begin{equation*}
c = \frac{n-1}{2(m(n-m)+1-n)}
\end{equation*}

Now we can calculate values of entries in $A_U$.  For example, consider the entry of $A_U$ whose  coordinates
correspond to variables $x_{ij}$ and $x_{ip}$ for $i,j \in U$ and $p \in
V-U$.   Recall that this entry is calculated in equation \eqref{entryA_U} for $I = \{\{i,j\},\{i,p\}\}$.  Note that 
\begin{equation*}
\sum_{\substack{u \in U \\ v \in V-U}}\sum_{\substack{x \in X \\ I \subset x}}x_{uv}
\end{equation*}
counts the number of Hamiltonian cycles containing the edges $\{i,j\},\{i,p\}$ and some edge from $U$ to $V-U$, and 
\begin{equation*}
\sum_{\substack{x \in X \\ I \subset x}}2
\end{equation*}
is simply 2 times the number of Hamiltonian cycles containing edges $\{i,j\}$ and $\{i,p\}$.  Using Lemma \ref{paths}, we can see
\begin{equation*}
\sum_{\substack{x \in X \\ I \subset x}}2 = 2(n-3)!
\end{equation*}
Note that $\{i,p\}$ is an edge from $U$ to $V-U$, and as was just stated, Lemma \ref{paths} tells us that the number of Hamiltonian cycles containing $I$ and edge $\{i,p\}$ (i.e. the number of Hamiltonian cycles containing $I$) is $(n-3)!$.  

There are $m-1$ edges from $U$ to $V-U$ containing vertex $p$ but not vertex $i$.  One of these edges, namely edge $\{j,p\}$, is not in any Hamiltonian cycles containing $I$ (because $n \not= 3$).  For each of the other edges, there are $(n-4)!$ Hamiltonian cycles containing that edge and $I$.

There are $n-m-1$ edges from $U$ to $V-U$ containing vertex $i$ but not vertex $p$.  However, none of these are in a Hamiltonian cycle containing $I$, because the vertex $i$ must have exactly 2 edges incident to it in a Hamiltonian cycle.

There are $n-m$ edges from $U$ to $V-U$ containing vertex $j$.  One of them (again, edge $\{j,p\}$) is not in any Hamiltonian cycle containing $I$.  For each of the other edges, there are $(n-4)!$ Hamiltonian cycles containing that edge and $I$.

And finally, there are $(m-2)(n-m-1)$ edges from $U$ to $V-U$ which do not contain any of vertices $i,j$, or $p$.  For each of these edges, Lemma \ref{paths} tells us that there are $2(m-4)!$ Hamiltonian cycles containing that edge and $I$.

Thus, from all of these arguments, we can calculate the value of $A_U$ whose coordinates correspond to variables $x_{ij}$ and $x_{ip}$:
\begin{align}
& \frac{c}{|X|}\sum_{\substack{u \in U \\ v \in V-U}}\sum_{\substack{x \in X  \\
 \{i,j\} \subset x \\ \{i,p\} \subset x}}x_{uv} +\sum_{\substack{x \in X \\ \{i,j\} \subset x \\ \{i,p\} \subset x }}2
\notag \\
 = &  c \frac{2}{(n-1)!} \bigg(\sum_{\substack{u \in U \\ v \in V-U}}\#\{\text{Hamiltonian paths
containing } 
 \{i,j\}, \{i,p\}, \{u,v\}\} \notag \\
& \quad \quad \quad \quad -2(n-2)!\bigg)  \notag \\
=& \frac{n-1}{2(m(n-m)+1-n)} \frac{2}{(n-1)!} \bigg( (n-3)! +(m-2)(n-4)!  \notag\\
& \quad \quad \quad \quad  +(n-m-1)(n-4)!+(m-2)(n-m-1)2(n-4)!-2(n-3)!\bigg) \notag \\
= & \frac{2(m-2)}{(n-2)(n-3)(m-1)} \label{example}
\end{align}  

All other entries of $A_U$ are found analogously.  

Now we turn back to the eigenvectors of $A_U$.  Suppose $m \geq 4$ and let $i,j,f, g$ be four distinct vertices in $U$. 
Consider the vector $u = x_{ij}-x_{jf}+x_{fg}-x_{g i}$.  Note that the coefficient of 1 in the vector $A_Uu$ is
\begin{multline}
\frac{c}{|X|}\bigg(\sum_{\substack{u \in U \\ v \in V-U}}\sum_{\substack{x \in X \\ \{i,j\} \subset x}}x_{uv} - \sum_{\substack{x \in X \\ \{i,j\} \subset x}}2 \bigg) - \frac{c}{|X|}\bigg(\sum_{\substack{u \in U \\ v \in V-U}}\sum_{\substack{x \in X \\ \{j,f\} \subset x}}x_{uv} - \sum_{\substack{x \in X \\ \{j,f\} \subset x}}2 \bigg) \\
+ \frac{c}{|X|}\bigg(\sum_{\substack{u \in U \\ v \in V-U}}\sum_{\substack{x \in X \\ \{f,g\} \subset x}}x_{uv} - \sum_{\substack{x \in X \\ \{f,g\} \subset x}}2 \bigg) - \frac{c}{|X|}\bigg(\sum_{\substack{u \in U \\ v \in V-U}}\sum_{\substack{x \in X \\ \{g,i\} \subset x}}x_{uv} - \sum_{\substack{x \in X \\ \{g,i\} \subset x}}2 \bigg) \label{const1}
\end{multline}
Note that the permutation $\pi = (if)$ permutes only vertices in $U$.  Thus, by an earlier remark, we know that the
entries of $A_U$ corresponding to $I = \{ij\}$ and $\{j,f\}$ are equal, as are the entries of $A_U$ corresponding to
$I = \{f,g\}$ and $\{i,g\}$.  By inspecting
\eqref{const1}, we can see that this implies that \eqref{const1} is equal to 0.

The coefficient of $x_{ab}$ in the vector $A_Uu$ is 
\begin{multline}
\frac{c}{|X|}\bigg(\sum_{\substack{u \in U \\ v \in V-U}}\sum_{\substack{x \in X \\ \{a,b\}, \\ \{i,j\} \subset x}}x_{uv} - \sum_{\substack{x \in X \\ \{a,b\}, \\ \{i,j\} \subset x}}2 \bigg) - \frac{c}{|X|}\bigg(\sum_{\substack{u \in U \\ v \in V-U}}\sum_{\substack{x \in X \\ \{a,b\}, \\ \{j,f\} \subset x}}x_{uv} - \sum_{\substack{x \in X \\ \{a,b\}, \\ \{j,f\} \subset x}}2 \bigg) \\
+ \frac{c}{|X|}\bigg(\sum_{\substack{u \in U \\ v \in V-U}}\sum_{\substack{x \in X \\ \{a,b\}, \\ \{f,g\} \subset x}}x_{uv} - \sum_{\substack{x \in X \\ \{a,b\}, \\ \{f,g\} \subset x}}2 \bigg) - \frac{c}{|X|}\bigg(\sum_{\substack{u \in U \\ v \in V-U}}\sum_{\substack{x \in X \\ \{a,b\}, \\ \{g,i\} \subset x}}x_{uv} - \sum_{\substack{x \in X \\ \{a,b\}, \\ \{g,i\} \subset x}}2 \bigg) \label{var1}
\end{multline}
Again, we note that for any $\pi \in S_m \times S_{n-m}$ ($\pi$ is a permutation of the vertices, leaving $U$ and $V-U$ fixed)
 we know that the entry of $A_U$ corresponding to edges in $I$ is equal to the entry of $A_U$ corresponding to edges
in $\pi I$, for any 0, 1, or 2 element set of edges $I$. 
Thus, upon inspection, we can see that this implies that
\eqref{var1} is 0 for any $\{a,b\}$ unless $\{a,b\}$
is one of $\{i,j\}, \{j,f\}, \{f,g\},$ or $\{g,i\}$.  Finally, we calculate that the coefficient of $x_{ij}$ is 
\begin{gather*}
\frac{c}{|X|}\bigg(\sum_{\substack{u \in U \\ v \in V-U}}\sum_{\substack{x \in X \\  \{i,j\} \subset x}}x_{uv} - \sum_{\substack{x \in X  \\ \{i,j\} \subset x}}2 \bigg) - \frac{c}{|X|}\bigg(\sum_{\substack{u \in U \\ v \in V-U}}\sum_{\substack{x \in X \\ \{i,j\}, \\ \{j,f\} \subset x}}x_{uv} - \sum_{\substack{x \in X \\ \{i,j\}, \\ \{j,f\} \subset x}}2 \bigg) \\
+ \frac{c}{|X|}\bigg(\sum_{\substack{u \in U \\ v \in V-U}}\sum_{\substack{x \in X \\ \{i,j\}, \\ \{f,g\} \subset x}}x_{uv} - \sum_{\substack{x \in X \\ \{i,j\}, \\ \{f,g\} \subset x}}2 \bigg) - \frac{c}{|X|}\bigg(\sum_{\substack{u \in U \\ v \in V-U}}\sum_{\substack{x \in X \\ \{i,j\}, \\ \{g,i\} \subset x}}x_{uv} - \sum_{\substack{x \in X \\ \{i,j\}, \\ \{g,i\} \subset x}}2 \bigg) \\
= \frac{n-1}{2(m(n-m)+1-n)} \frac{2}{(n-1)!} \bigg( 2(n-m)(n-3)!+(m-2)(n-m)2(n-3)! \\
-2(n-2)!\bigg) \\
- \frac{n-1}{2(m(n-m)+1-n)} \frac{2}{(n-1)!} \bigg( 2(n-m)(n-4)!+(m-3)(n-m)2(n-4)! \\
-2(n-3)!\bigg) \\
 + \frac{n-1}{2(m(n-m)+1-n)} \frac{2}{(n-1)!} \bigg( 4(n-m)2(n-4)!+(m-4)(n-m)4(n-4)! \\
 -2\cdot2(n-3)!\bigg) \\
 - \frac{n-1}{2(m(n-m)+1-n)} \frac{2}{(n-1)!} \bigg( 2(n-m)(n-4)!+(m-3)(n-m)2(n-4)! \\
 -2(n-3)!\bigg) \\
 =\frac{2(m-2)}{(n-2)(m-1)}
\end{gather*}
(where, again, the above calculations for entries $A_U$ are analogous to the calculation in \eqref{example}).  Similarly, the entry of $A_Uu$ corresponding to $x_{jf}$ is $-\frac{2(m-2)}{(n-2)(m-1)}$, corresponding to $x_{fg}$ is $\frac{2(m-2)}{(n-2)(m-1)}$, and corresponding to $x_{g i}$ is $-\frac{2(m-2)}{(n-2)(m-1)}$.  

Thus, we have proven the first two rows in our eigenvalue table for Lemma \ref{h_U}.  The remaining eigenvalues are
proven analogously, by simply inspecting the eigenvectors, and using Lemma \ref{paths} to analyze how $A_U$ acts on
the eigenvectors. 

\end{proof}

We note that in the case $m=2$, $U$ consists of two vertices, say $U = \{i,j\}$.  In this case, if the polynomial $p(x) = x_{ij}$ is
 nonzero then the function $h_U$ is 0.  Thus, the row and column of $A_U$ corresponding to $x_{ij}$ will be 0, giving
another  eigenvector corresponding to eigenvalue 0.  Hence, if $|U|=2$, then $h_U$ already lies on the boundary of
$P_1$

\begin{proof}[Proof of Lemma \ref{ones}]
Just as with Lemma \ref{h_U}, we simply calculate the action of $A_{\mathbbm{1}}$ on the eigenvectors.  Because it is 
completely analogous, this calculation is omitted.

\end{proof}

\begin{proof}[Proof of Theorem\ref{eigenvalues}]
Here we need to calculate the eigenvectors for the matrix $A_{U,a}=\alpha A_U + (1-a)A_\mathbbm{1}$
corresponding to the quadratic form $q_{ah_U+(1-a)\mathbbm{1}}$.  Using Lemmas
\ref{h_U} and
\ref{ones}, we can see that each eigenvector that we found for $A_U$ is also an eigenvector for $A_\mathbbm{1}$.  For
each of eigenvalues that we found, when we calculate the dimension of the span of the associated eigenvectors, we see
that we already know the following eigenvalues and multiplicities for $q_{ah_U+(1-a)\mathbbm{1}}$:

\begin{tabular}{ll}
\textbf{Multiplicity} & \textbf{Eigenvalue} \\
n & 0 \\
$\frac{m(m-3)}{2}$ &  $a\frac{2(m-2)}{(n-2)(m-1)} + (1-a)\frac{2}{n-1}$ \\
$\frac{(n-m)(n-m-3)}{2}$ &  $a\frac{2(n-m-2)}{(n-2)(n-m-1)} + (1-a)\frac{2}{n-1}$ \\
$(n-m-1)(m-1)$ & $a\frac{2(mn^2-nm^2-n^2+4n-3mn+3m^2-4)}{(n-2)(n-3)(mn-m^2-n+1)} + (1-a)\frac{2}{n-1}$ \\
$m-1$ & $a\frac{2(m-2)}{(n-3)(m-1)} + (1-a)\frac{2}{n-1}$ \\
$n-m-1$& $a\frac{2(n-m-2)}{(n-3)(n-m-1)} + (1-a)\frac{2}{n-1}$
\end{tabular}

Thus, the total number of eigenvalues (with multiplicities) that we know so far is
\begin{gather*}
n+\frac{m(m-3)}{2}+\frac{(n-m)(n-m-3)}{2}+(n-m-1)(m-1) \\
+(m-1)+(n-m-1) 
=\frac{n(n-1)}{2}-1
\end{gather*}
Since the dimension of the space of
linear functions on $\R^{n(n-1)/2}$ is $\frac{n(n-1)}{2} +1$, there are two eigenvalues yet to calculate.  From
Theorem \ref{trace}, we can already calculate the trace of $aA_U+(1-a)A_\mathbbm{1}$.  Using the same techniques as in
the proof of Lemma \ref{h_U}, we can calculate all of the entries of $aA_U+(1-a)A_\mathbbm{1}$, and use these to
calculate the diagonal entries, and thus the trace of $(aA_U+(1-a)A_\mathbbm{1})^2$.  Using this information, the
remaining two eigenvalues that we find are
\begin{equation*}
\frac{c+\sqrt{d}}{2(mn^3-n^3-5mn^2+6n^2-m^2n^2-11n+5m^2n+6mn+6-6m^2)(n-1)}
\end{equation*}
and
\begin{equation*}
\frac{c-\sqrt{d}}{2(mn^3-n^3-5mn^2+6n^2-m^2n^2-11n+5m^2n+6mn+6-6m^2)(n-1)}
\end{equation*}
where
\begin{multline*}
c =3n^4m-3n^4-2an^3-3m^2n^3+17n^3-14mn^3+4an^2+8amn^2-27n^2+14m^2n^2 \\+13mn^2-13m^2n-16amn-8am^2n+6mn+2an+7n-6m^2+16am^2-4a+6
\end{multline*}
and
\begin{align*}
d=&(2-3n+n^2)(162-72a-675n+324mn-324m^2-120an^3-12an^5\\
&+72an^4 
+333m^2n^3+702m^3n^2-288am^4+252am^2n^3+9n^6m^2+8a^2\\
& -4a^2n^3+4a^2n^4
+200a^2m^4 -12a^2n^2+4a^2n-136a^2m^2-104a^2m^2n^3\\
& +80a^2m^2n^2-400a^2m^3n  +136a^2mn 
+256a^2m^3n^2+24a^2m^2n^4\\
&+232a^2m^2n - 24a^2mn^4-232a^2mn^2+120a^2mn^3-128a^2m^4n 
\\
& +24a^2m^4n^2 -48a^2m^3n^3  -60am^2n^2+576am^3n-360amn-351m^4n\\
& +180n^5m-18n^6m 
 +432n^4 -684n^4m-954n^3+162m^4-324m^3n\\
 & -63m^2n^5-480am^3n^2+12amn^5 +1125n^2 
+132an+360am^2-522m^3n^3\\
& +261m^4n^2-99n^5 +9n^6-1026mn^2-1062m^2n^2+1224mn^3 
+1026m^2n\\
& -60am^2n^4-588am^2n-12amn^4+588amn^2-228amn^3+240am^4n 
\\
& -48am^4n^2+96am^3n^3-18m^3n^5+9n^4m^4 +81n^4m^2+162n^4m^3\\
 &-81n^3m^4)
\end{align*}
We note that if $a = \sqrt{n}$, then for any $m \leq \frac{n}{2}$, the second eigenvalue listed above is $\leq 0$.

\end{proof}

\section*{Acknowledgements}
The author would like to thank Alexander Barvinok for his helpful ideas and encouragement, see section \ref{cone1}.

\bibliographystyle{plain}
\bibliography{semidefiniteTSP}
\end{document}